\documentclass[10pt]{amsart}

\usepackage{epsfig,amsmath,amsthm, amssymb}

\newtheorem{theorem}{Theorem}[section]
\newtheorem{proposition}[theorem]{Proposition}

\newtheorem{corollary}[theorem]{Corollary}

\numberwithin{equation}{section}

\title{Quotients of hypersurfaces in weighted projective space}
\author{GILBERTO BINI}

\email{gilberto.bini@unimi.it} \curraddr{{\sc Dipartimento di
Matematica, Universit\`a degli Studi di Milano  \\ Via C. Saldini 50 -
20133 Milano  Italy}}


\begin{document}

\maketitle

\begin{abstract} In \cite{bvg} some quotients of one-parameter families of Calabi-Yau varieties are related to the family of Mirror
Quintics by using a construction due to Shioda. In this paper, we
generalize this construction to a wider class of varieties. More
specifically,  let $A$ be an invertible matrix with non-negative
integer entries. We introduce varieties $X_A$ and $\overline{M}_A$
in weighted projective space and in ${\mathbb P}^n$, respectively.
The variety $\overline{M}_A$ turns out to be a quotient of a Fermat
variety by a finite group. As a by-product, $X_A$ is a quotient of a
Fermat variety and $\overline{M}_A$ is a quotient of $X_A$ by a
finite group. We apply this construction to some families of
Calabi-Yau manifolds in order to show their birationality.
\end{abstract}

\section{Introduction}\label{sec1}

Hypersurfaces in weighted projective space have been investigated by
many authors in connection with Physics, in particular Mirror
Symmetry: see, for instance, \cite{cdk}, \cite{cls} and
\cite{Yau}. The usual quintic threefold is an example of
Calabi-Yau manifold in ordinary projective space. Its mirror can
be described in terms of quotients of a one-parameter family of
quintics, the Dwork pencil. In \cite{bvg}, we investigated other
one-parameter families of Calabi-Yau manifolds and related them to
the family of Mirror quintics. Our main tool was a construction due
to Shioda \cite{shio}, which clarified previous work in \cite{dgj}.

In the present paper we generalize this construction to
hypersurfaces in weighted projective space. Let $A$ be an invertible
matrix of size $n$ with non-negative integer entries. For such a
matrix we define a weighted homogeneous polynomial $F_A$, a sum of
$n$ monomials. Its zero locus gives a projective scheme $X_A$ in
weighted projective space $W{\mathbb P}^{n-1}(q_1, \ldots, q_n)$.
The weights are determined by the relation $A{\bf q}=d {\bf e}$,
where ${\bf q}=(q_1, \ldots, q_n)$, ${\bf e }=(1,1, \ldots, 1)$ and
$d$ is the smallest positive integer such that $dA^{-1}$ has integer
entries. If the total degree of $F_A$ equals the degree of the
anticanonical bundle and the singularities of $X_A$ are canonical,
then $X_A$ is a Calabi-Yau manifold as explained in Section
\ref{sec2}. It is therefore interesting to study the quotients of
$X_A$ in relation with the mirror varieties of $X_A$. For this
purpose, we introduce a manifold $\overline{M}_A \subset {\mathbb
P}^n$, which we refer to as the Shioda quotient. It is the image of
a suitable map (first introduced by Shioda in \cite{shio} and
applied to a similar context in \cite{bvg}) of a Fermat variety.
Actually, we prove that the Shioda quotient is the quotient of a
Fermat variety by a finite group. As a by-product, it is also the
quotient of $X_A$ by a finite group. We describe these groups in
detail.

It is natural to investigate deformation families of $X_A$ defined
by $F_{A,t}=F_A-tx_1\ldots x_n$, where  $x_i$ are variables of
degree $q_i$. They are particularly meaningful when the Hodge number
$h^{2,1}$ is one, since then $F_{A,t}$ gives a versal deformation. Some examples of deformation families of $X_A$
when $X_A$ is a Calabi-Yau manifold in weighted projective space are
given in \cite{sch}. Their mirror all have $h^{2,1}=1$. We investigate the birational classes associated to these families and prove they are birational to the mirror families of the first
four hypergeometric Calabi-Yau families studied by Rodriguez Villegas \cite{RV}
and listed in \cite{mey}, p. 134.

\section{Preliminaries}\label{sec2}

\subsection{Weighted projective spaces} We briefly recall the definition of weighted projective space. Let $(q_1, \ldots, q_n)$ be a sequence of positive integers. As customary, set
$$
W{\mathbb P}^{n-1}(q_1, \ldots, q_n):= ({\mathbb C}^{n} \setminus \{0\})/\sim,
$$
where the equivalence relation is
$$
(x_1, \ldots, x_n) \sim (\lambda^{q_1}x_1, \ldots, \lambda^{q_n}x_n)
$$
whenever $\lambda \in {\mathbb C} \setminus \{0\}$. Recall that
$$
W{\mathbb P}^{n-1}(q_1, \ldots, q_n) \cong W{\mathbb
P}^{n-1}(q_1/a_1, \ldots, q_n/a_n),
$$
where $a_i=l.c.m.(d_1, \ldots, \widehat{d_i}, \ldots, d_n)$ and
$d_i=g.c.d.(q_1, \ldots, \widehat{q_i}, \ldots, q_n)$. Weighted
projective space are almost always singular. As proved in \cite{Di},
the singular locus of $W{\mathbb P}^{n-1}(q_1, \ldots, q_n)$ can be
described in the following way. Let $p$ be a prime. Let $I(x)=\{j;x_j \neq 0\}$. Then define
$$
Sing_p(W{\mathbb P}^{n-1}):=\{ x \in W{\mathbb P}^{n-1}(q_1,\ldots,
q_n) : p|q_i \, \hbox{for any}\, i \in I(x)\}.
$$

The singular locus of the weighted projective space is given by the
union over all primes of $Sing_p(W{\mathbb P}^{n-1})$.

 As explained, for instance, in \cite{CK},  weighted projective space is a toric variety.  Set $Q:=\sum_i q_i$. We recall that the canonical sheaf of $W{\mathbb P}^{n-1}(q_1, \ldots, q_n)$ is given by ${\mathcal O}( - Q)$, which is not always a line bundle. As proved in \cite{CK}, Lemma 3.5.6., the canonical sheaf is a line bundle if and only if $q_i |Q$ for all $i=1, \ldots, n$. Under this assumption, weighted projective space is a Fano toric variety.

Some of our hypersurfaces are Calabi-Yau varieties. Following \cite{CK}, a (possibly singular) Calabi-Yau variety is an $m$-dimensional normal compact variety $X$ which satisfies the following conditions: \\
(i) $X$ has at most Gorenstein canonical singularities; \\
(ii) the dualizing sheaf of $X$ is trivial;\\
(iii) $H^1(X, {\mathcal O}_X)=\ldots = H^{m-1}(X, {\mathcal O}_X)=0$.

  Let $f$ be a weighted homogeneous polynomial such that the zero locus $\{f=0\}$ is quasi-smooth (according to Definition 3.1.5 in \cite{Dol}).
  Since a quasi-smooth scheme has finite quotient singularities (\cite{Dol}, Thm. 3.1.6), the locus
   $\{f=0\}$ is normal and Cohen-Macaulay; furthermore, it is Gorenstein with dualizing sheaf ${\mathcal O}(d)$,
   where $d$ is the weighted degree of $f$. Assume that $\{f=0\}$ has at most canonical singularities. Following the arguments in
   Proposition 4.1.3 in \cite{CK}, it follows that $\{f=0\}$ is a Calabi-Yau variety when $d=Q$. Notice that if there exists a crepant
   resolution of $\{f=0\}$, the singularities are canonical by
   definition.

\subsection{The Shioda maps} Let $A$ be an invertible matrix with non-negative integer entries:
$$
A=(a_{ij}) \quad (\in M_{n}({\mathbb Z})), \quad a_{ij} \in {\mathbb Z}_{\geq 0}, \quad det(A) \neq 0.
$$

For such a matrix we define a polynomial in $n$ variables, a sum of $n$ monomials:
$$
F_A:=\sum_{i=1}^n\prod_{j=1}^n x_j^{a_{ij}}= x_1^{a_{11}}x_2^{a_{12}}\ldots x_n^{a_{1n}} + x_1^{a_{21}}x_2^{a_{22}}\ldots x_n^{a_{2n}} + \ldots
$$

Let $d$ be the smallest positive integer such that $B:=dA^{-1}$ is in $M_n({\mathbb Z})$. Then set
\begin{equation}
\label{qu}
{\bf q}:= B{\bf e},
\end{equation}
where ${\bf e}=(1, \ldots, 1)$ and ${\bf q}=(q_1, \ldots, q_n)$. Clearly, this implies that
\begin{equation}
\label{di}
A{\bf q}=d{\bf e}.
\end{equation}

The zero locus $X_A=Z(F_A)$  is a (not necessarily smooth or irreducible) projective variety $X_A$, which is contained
in $W{\mathbb P}^{n-1}(q_1, \ldots, q_n)$. Let $m$ be the greatest common divisor of the $q_i$'s. Define $a_i=q_i/m$.
Thus we have $$X_A \subset W{\mathbb P}^{n-1}(q_1, \ldots, q_n) \cong W{\mathbb P}^{n-1}(a_1, \ldots, a_n).$$

We have a rational map $\phi_A$ from ${\mathbb P}^{n-1}$ of degree $d$ to $X_A$ defined by:
$$
\phi_A: {\mathbb P}^{n-1} \rightarrow W{\mathbb P}^{n-1}(q_1, \ldots, q_n),$$
  $$(y_1: \ldots : y_n) \rightarrow (x_1: \ldots : x_n), \quad x_j=\prod_{k=1}^ny_k^{b_{jk}}.
 $$

Notice that each $y_j$ has degree one, so $deg(x_j) = \sum_k b_{jk}=q_j$. Hence $\phi_A$ is indeed a rational map from ${\mathbb P}^{n-1}$
to $W{\mathbb P}^{n-1}(q_1, \ldots, q_n)$.

Assume further that $d=Q=\sum_jq_j$, so that the variety $X_A$ gives
a Calabi-Yau provided the singularities are canonical.  We can read
this condition on the coefficients of the matrix $A^{-1}$. In fact,
we have:
$$
d= Q= \sum_j q_j ={}^t{\bf e}{\bf q}= d ^t{\bf e}A^{-1}{\bf e},
$$
which gives
\begin{equation}
\label{boh}
^t{\bf e}A^{-1}{\bf e}=1;
\end{equation} in other words, $\sum_{ij}a_{ij}'=1$, where $a'_{ij}$ are the entries of $A^{-1}$.

We define a rational map
$$
q_A:W{\mathbb P}^{n-1}(q_1, \ldots, q_n) \rightarrow \overline{Im(q_A)}:=\overline{M}_A \subset {\mathbb P}^n
$$
in the following way:
$$
(x_1:\ldots :x_n) \rightarrow (u_0: u_1 : \ldots : u_n):=\left(\prod _{j=1}^n x_j : \prod_{j=1}^nx_j^{a_{1j}}:
 \ldots : \prod_{j=1}^n x_j^{a_{nj}}\right).
$$

Since $deg(x_i)=q_j$ and $A{\bf q}= d{\bf e}$, we have
$$
deg(u_0)=\sum_j q_j=d, \qquad deg(u_k)=\sum_j a_{kj} q_j =d,
$$
hence $q_A$ is well-defined.

Finally, we describe the composition $q_A \circ \phi_A$, which will be used in the next section. First, as $BA=AB=dI_n$ (where $I_n$ is the identity matrix) we have $\sum_j a_{lj}b_{jk}=d\delta_{lk}$, where $\delta_{kl}$ is the Kronecker delta.
 Second, we set
$$
{\bf q}':=d\,{} ^tA^{-1}{\bf e}= d \, ^t {\bf e} A^{-1} = \, ^t {\bf e} B;
$$
so $q_k'=\sum_j b_{jk}.$ This said, it is easy to check that the composition $q_A \circ \phi_A: X_d \subset {\mathbb P}^{n-1} \rightarrow {\overline M}_A \subset {\mathbb P}^n$ is given by
\begin{equation}
\label{composition} (u_0: u_1 : \ldots : u_n) = \left(
\prod_{k=1}^ny_k^{q'_k}: y_1^d : \ldots : y_n^d \right),
\end{equation}
where $X_d$ is the Fermat variety $\{ \sum_{i=1}^n y_i^d=0\}$ and
${\overline M}_A$ is the closure of the image of $X_d$.

Let us consider the projection
\begin{equation}
\label{forget}
\begin{array}{lcr}
\pi: {\overline M}_A \subset {\mathbb P}^n & \rightarrow & V \cong {\mathbb P}^n \\
(u_0: u_1: \ldots : u_n) & \rightarrow & (u_1 : \ldots : u_n),
\end{array}
\end{equation}
where $V$ is the closure of $\pi({\overline M}_A)$.  It is easy to check
that $V$ is isomorphic to the ${\mathbb P}^{n-1}$ given by
$\sum_{i=1}^n u_i=0$ as the Fermat equation has to be satisfied.

From now on, we will assume that $q_k'$  is strictly positive for
any $k$. By direct inspection, we obtain the following equations for
the image of the Fermat variety $X_d$ under $q_A \circ \phi_A$:
\begin{equation}
u_1 + \ldots + u_n =0, \quad u_0^d= u_1^{q_1'} \ldots u_n^{q_n'}.
\label{first}
\end{equation}

Set $m':=g.c.d(d, q_1', \ldots, q_n')$. Hence $d=m'a'$ and
$q_k'=a_k'm'$, so the composition $q_A \circ \phi_A$ is given by
\begin{equation}
\label{newim}
(u_0:\ldots:u_n)=\left(\prod_{k=1}^ny_k^{q'_k}: y_1^d : \ldots : y_n^d \right)=\left(
\prod_{k=1}^n(y_k^{m'})^{a'_k}: (y_1^{m'})^{a'} : \ldots : (y_n^{m'})^{a'} \right).
\end{equation}

By composing \eqref{newim} with the map $t_k=y^{m'}_k$ for $k=1, \ldots, n$, we get
$$
(u_0:\ldots:u_n)=\left(\prod_{k=1}^nt_k^{a'_k}: t_1^{a'} : \ldots : t_n^{a'} \right),
$$
so the equations defining ${\overline M}_A$ are the following:
\begin{equation}
\label{second}
u_1+ \ldots + u_n=0, \quad u_0^{a'}= u_1^{a_1'} \ldots u_n^{a_n'}.
\end{equation}

 In the next section, we will show that under our assumptions on the $q_j$'s, the equation \eqref{second} define a very singular  model in $(n-1)$-dimensional projective space for  a manifold with $h^{n-2,0}=1$ of degree $a' ( > n$ in general).

\section{The Shioda quotient}\label{sec3}

In this section we assume that $A \in M_n({\mathbb Z}_{\geq 0})$ is an invertible matrix such that \eqref{boh} holds. We will introduce "natural" automorphism groups and study the quotients by these groups.

\subsection{The automorphism groups}
 Let $\zeta= \zeta_d$ be a generator of the cyclic group of $d$-th roots of unity, where $d$ is the smallest positive integer such that $dA^{-1}$ has integer entries.
For ${\bf k}=(k_1, \ldots, k_n) \in ({\mathbb Z}/d{\mathbb Z})^n$ we define an automorphism $g_{{\bf k}}$ of ${\mathbb P}^{n-1}$ by
$$
g_{{\bf k}}(y_1:\ldots : y_n):=\left(\zeta^{k_1}y_1 : \ldots : \zeta^{k_n} y_n\right).
$$
Note that ${\bf a},{\bf b} \in \mu_d^n$ define the same automorphism iff ${\bf a} - {\bf b} =
(k, \ldots, k)$ for some $k \in \mu_d$. Define $\Gamma_d$ to be the
quotient group
$$
\Gamma_d:=\mu_d^n/\langle g_{(1,1, \ldots, 1)} \rangle  \, \, \left ( \subset Aut({\mathbb P}^{n-1})\right)
$$

Notice that $\Gamma_d \cong \mu_d^{n-1}$, hence $\# \Gamma_d = d^{n-1}$. The group $\Gamma_d$ is a
subgroup of the automorphism group of the Fermat variety $X_d$.

The map $q_A \circ \phi_A$ \eqref{composition} is invariant under
the subgroup of $\Gamma_d$ given by
$$
\Gamma({\bf q'}):=\left\{ g_{{\bf k}}: {\bf k}= (k_1, \ldots,
k_n); \sum_j k_j q'_j \equiv 0 \, \, \hbox{mod} \, \, {d}\right\}\Biggl/ \langle g_{(1,\ldots, 1)} \rangle.
$$

In other words, we have
$$
(q_A \circ \phi_A)(g_{{\bf k}}(y_1: \ldots: y_n))=(q_A \circ \phi_A)
(y_1: \ldots: y_n)
$$
for all $g_{{\bf k}}$ in $\Gamma({\bf q}')$ and all $(y_1, \ldots, y_n) \in X_d$.
Notice that $g_{(1, \ldots, 1)}$ is an element of $\Gamma({\bf q'})$. In fact, by \eqref{boh}:
$$
\sum_j q'_j= \, ^t{\bf e}{\bf q}'=d\, ^t{\bf e}\, ^tA^{-1}{\bf e}=d.
$$

The coordinate functions of the Shioda map
$\phi_A$ are products of the $y_i$. If $\phi_A(y_1: \ldots: y_n)=(x_1: \ldots : x_n)$, then
$$
\phi_A(g_{{\bf k}}(y_1: \ldots: y_n))= \left( \zeta^{k'_1}x_1:
\ldots : \zeta^{k'_n}x_n \right).
$$

As $x_j= \prod y_k^{b_{jk}}$, the column vector ${\bf k}'$ is
obtained from the column vector ${\bf k}$ as ${\bf k}'=B{\bf
k}$. Thus we get a homomorphism
$$
\Gamma({\bf q}') \rightarrow Aut(X_A), \quad g_{\bf k} \rightarrow g_{B{\bf k}},
$$
which is well defined since $B{\bf e} \equiv 0$ mod $d$.

The kernel (image resp.) of this homomorphism will be denoted by $\Gamma_A$ ($H_A$ resp.).
Notice that $\Gamma_A$ is the subgroup of $\Gamma({\bf q}')$, which is generated by the images of the $g_{{\bf k}}$
such that $B{\bf k} \equiv 0$ mod $d$.

\subsection{The birational model}   We recall that two rational maps between
algebraic varieties $f_i: X \rightarrow Y_i$ for $i=1,2$ are said to be birationally
 equivalent if there is a Zariski open subset $U$ of $X$ and
 there are Zariski open subsets $U_i \subset Y_i$ with
 an isomorphism $\phi: U_1 \rightarrow U_2$ such that $\phi \circ f_1 = f_2$ on $U$.

\begin{theorem}
\label{central}
Let $A$ be an invertible $n \times n$ matrix with integer entries such that $X_A$ is irreducible and \eqref{boh} holds. Then the
composition $q_A \circ \phi_A$ is birational to the quotient map $X_d \rightarrow X_d/\Gamma({\bf q}')$; hence $X_d/\Gamma({\bf q'})$ is birational to ${\overline M}_A$.
\end{theorem}
\proof

The composition $q_A \circ \phi_A$ is given by \eqref{composition}.
Also, recall the map $\pi$ defined in  \eqref{forget}. The
composition of $ q_A \circ \phi_A $ and $\pi$ yields a map from the
Fermat variety $X_d$ to $V \cong {\mathbb P}^{n-1}$, which
corresponds to an abelian extension  with group $\Gamma_d$ of
function fields - recall that $X_d$ is the Fermat variety and $u_i =
y_i^d$. This means that $X_d \rightarrow V$ is the quotient for the
group $\Gamma_d$, namely $X_d/\Gamma_d = V$. Therefore, by abelian
Galois theory, each subfield is obtained as an invariant field under
a finite subgroup of $\Gamma_d$. Thus, the map $X_d \rightarrow
\overline{M}_A$ corresponds to a quotient by a finite subgroup of
$\Gamma_d$. Now, we show that this subgroup is isomorphic to
$\Gamma({\bf q}')$. The map $\pi$ corresponds to an abelian
extension of function fields with group ${\mathbb Z}/a'{\mathbb Z}$,
where $d=m'a'$ and $m'=g.c.d.(d, q_1', \ldots, q_n')$.  The group
$\Gamma_d$ acts on ${\overline M}_A$ only through the variable
$u_0$, and $g_{\bf k}:u_0 \rightarrow \zeta_{m'}^{{\bf q}'\cdot {\bf
k}}u_0$, where $\zeta_{m'}=\zeta^{a'}$ is a primitive $m'$-th root
of unity. The kernel of this action is exactly $\Gamma({\bf q'})$,
hence the map $X_d \rightarrow {\overline M}_A$ corresponds to an
extension of function fields with group $\Gamma({\bf q'})$. Hence,
the claim follows.
\endproof

By the Theorem above, the order of the group $\Gamma({\bf q}')$ is $d^{n-2}m'$.

\begin{corollary}
The maps $\phi_A : X_d \rightarrow X_A$ and $q_A: X_A \rightarrow \overline{M}_A$ are birational to quotient maps. In particular, $X_A$ (${\overline M}_A$ resp.) is birational to $X_d/\Gamma_A$ ($X_A/H_A$ resp.).
\end{corollary}

\proof By Theorem \ref{central}, the composition of $\phi_A$ and $q_A$ is a quotient map, namely:
$$
X_d \rightarrow X_A \rightarrow {\overline M}_A \approx X_d/\Gamma({\bf q}').
$$

The proof follows easily from arguments similar to those in Theorem 2.6 in \cite{bvg}.
\endproof

Now, we prove that the equations \eqref{first} and \eqref{second} give a very singular model for a manifold with $h^{n-2,0}=1$.
Assume, $q_i' > 0$. We recall that the vector space of holomorphic $(n-2)$-forms on a smooth
hypersurface $X=Z(F)$ of degree $d \geq n$ in ${\mathbb P}^{n-1}$ has a basis of the form
$$
\omega_{{\bf b}, F}=Res_X \left(y_1^{b_1}\ldots y_n^{b_n} \frac{\sum_i ^n (-1)^i y_i dy_1 \wedge \ldots \wedge \widehat{dy_i}
\wedge \ldots \wedge dy_n}{F} \right),
$$
where ${\bf b}=(b_1, \ldots, b_n)$ and $b_i \in {\mathbb Z}_{\geq 0}$, $\sum_i b_i = d-n$.

\begin{proposition}
There exists a unique holomorphic $(n-2)$ form on any resolution of ${\overline M}_A$.
\end{proposition}

\proof The action of an element $g_{{\bf k}}$ in $\Gamma({\bf q'})$
on $\omega_{{\bf b}, F}$ is given by $\zeta^{\sum_i (b_i+1)k_i}$,
where $b_i \in {\mathbb Z}_{\geq 0}$. A form is invariant with
respect to $\Gamma({\bf q'})$ if and only if $b_i+1=q'_i$ for all
$i$. The unique invariant form $\omega_{{\bf b}, F}$ with vector
${\bf b}={\bf q'}- {\bf e}$ descends to a form on the quotient
$X_d/\Gamma({\bf q'}) \approx {\overline M}_A$.
\endproof

\subsection{Some examples}\label{sec7}

{\it Example A.}  Let us consider
$F_A:=x_1^5+x_2^{10}+x_3^{10}+x_4^{10}+x_5^2=0$  in weighted
projective space $W{\mathbb P}^4(2,1,1,1,5)$. It is easy to check
that the corresponding hypersurface is smooth and does not intersect
the singularities of $W{\mathbb P}^4(2,1,1,1,5)$, which are two
isolated points.

The  matrix $A$ is given by $diag(5,10,10,10,2)$. The matrix $B$ is $diag(2,1,1,1,5)$ and $d=10$. The condition ${}^t{\bf
e}A^{-1}{\bf e}=1$ is satisfied. Moreover,
$$
{\bf q}=(2,1,1,1,5), \quad {\bf q}'=(2,1,1,1,5)
$$
so $q_i' >0$ for any $i=1, \ldots, 5.$ The equations cutting out
$\overline{M}_A$ in ${\mathbb P}^5$ are given by
$$
u_1+u_2+u_3+u_4+u_5=0,
$$
$$
u_0^{10}=u_1^2u_2u_3u_4u_5^5.
$$

The integer $d = 10$. Generators for the groups $\Gamma({\bf q}')$, $\Gamma_A$ and $H_A$ are as follows. Consider the elements ${\bf v_1}, \ldots , {\bf v_4} \in ({\mathbb Z}/10{\mathbb Z})^5$ given by:
$$
{\bf v}_1=(0,0,0,5,1), \qquad {\bf v}_2=(0,0,1,4,1),
$$
$$
{\bf v}_3=(1,0,0,3,1), \qquad {\bf v}_4=(0,1,0,4,1),
$$
then $^t {\bf q'}\cdot {\bf v}_i \equiv 0$ mod $10$ for $i=1, \ldots, ,4$. Moreover, we have ${\bf e}= 8{\bf v}_1+{\bf v}_2+{\bf v}_3+{\bf v}_4$; hence
$$
\Gamma({\bf q}') = \langle {\bf v}_1, {\bf v}_2, {\bf v}_3 \rangle \cong \mu_{10}^3.
$$

The group $\Gamma_A$ is isomorphic to $\mu_{10}$ and is generated by $(5,0,0,0,6)$. Finally, set
$$
{\bf w}_1=(0,0,1,4,5), \quad {\bf w}_2=(2,0,0,3,5), \quad {\bf w}_3=(0,1,0,4,5).
$$

As ${\bf w}_1 + {\bf w }_2 + {\bf w}_3 = {\bf q}$, we have
$$
H_A= \langle {\bf w}_1, {\bf w}_2 \rangle \cong \mu_{10}^2 \subset
Aut(X_A).
$$

The isomorphism $H_A \cong \mu_{10}^2$ was first suggested in
\cite{gpr}.  Finally, as mentioned in \cite{mor}, notice that
$h^{2,1}(X_A)=1$.

{\it Example B.} Let us consider the equation
$F_A:=x_1^{15}x_5+x_2^5+x_3^5+x_3x_4^5+x_2x_4^2=0$ in weighted
projective space $W{\mathbb P}^4(1,5,5,4,10)$. In this case, we do
not have a Fano toric veariety since $4$ does not divide
$25=1+5+5+4+10$. The zero locus does not intersect the singularities
of $W{\mathbb P}^4(1,5,5,4,10)$ and $F_A$ is a smooth variety. The
matrices $A$ and $B$ are given by
$$
A:=\left(
\begin{array}{ccccc}
15 & 0 & 0 & 0 & 1 \\
0 & 5 & 0 & 0 & 0 \\
0 & 0 & 5 & 0 & 0 \\
0 & 0 & 1 & 5 & 0 \\
0 & 1 & 0 & 0 & 2
\end{array}
\right),
 \quad B:= \left(
\begin{array}{ccccc}
10 & 1 & 0 & 0 & -5 \\
0 & 30 & 0 & 0 & 0 \\
0 & 0 & 30 & 0 & 0 \\
0 & 0 & -6 & 30 & 0 \\
0 & -15 & 0 & 0 & 75
\end{array}
\right).
$$

The integer $d$ equals $150$. Moreover, we have
$$
{\bf q}=6(1, 5, 5, 4, 10), \qquad {\bf q}'=(10,16,24,30,70)
$$

and $a'=g.c.d.(d, q_1', \ldots, q_5')=2$. Notice that $q_i' >0$ for any $i=1, \ldots, 5$.

A birational model of ${\overline M}_A$ is cut out by the equations
\begin{equation}
\label{notrans}
u_0^{75}=u_1^5u_2^8u_3^{12}u_4^{15}u_5^{35}, \quad u_1 + \ldots + u_5=0.
\end{equation}

Consider the vectors ${\bf r_i} \in {\mathbb Z}/d{\mathbb Z}$ given by:
$$
{\bf r}_1= (0,0,75,0,0), \qquad {\bf r}_2=(0,1,1,0,8),
$$
$$
{\bf r}_3=(1,0,0,0,2), \quad {\bf r}_4=(0,0,0,1,6), \quad {\bf r}_5=(0,0,5,0,9);
$$
then $9{\bf e}={\bf r}_1+9 {\bf r}_2+ 9{\bf r}_3 + 9{\bf r}_4- 15 {\bf r}_5$. The vectors ${\bf r}_i$ generate $\Gamma({\bf q}')$; in fact, the following holds:
$$
\Gamma({\bf q}')\cong \langle {\bf r}_1, {\bf r}_2, {\bf r}_3, {\bf r}_4, {\bf r}_5 \rangle/ < {\bf r}_1+9 {\bf r}_2+ 9{\bf r}_3 + 9{\bf r}_4- 15 {\bf r}_5 >\cong {\mathbb Z}/2{\mathbb Z} \times \left( {\mathbb Z}/150{\mathbb Z}\right)^3
$$

By using Magma, it is possible to check that
$$
\Gamma_A \cong \left( {\mathbb Z}/150{\mathbb Z}\right)^3,
$$
$$
H_A \cong {\mathbb Z}/2{\mathbb Z}.
$$

The group $H_A$ is generated by $g_{(0,75,75,0,75)}$, which maps
$(x_1: \ldots : x_5)$ to $(x_1: - x_2: - x_3 : x_4: -x_5)$. Recall
that $X_A \rightarrow {\overline M}_A \approx X_A/H_A$ is a double
cover.

 {\it Example C}. Consider the weighted homogeneous polynomial  $${F_A}:=x_1^2+x_2^3+x_3^{18}+x_4^{18}+x_5^{18}$$. It gives a quasi-smooth locus in weighted projective space $W{\mathbb P}^4(9,6,1,1,1)$.
 The matrices $A$ and $B$ are given by $diag(2,3,18,18,18)$ and
 $diag(9,6,1,1,1)$, respectively.

The integer $d$ is equal to $18$ and we have
$$
{\bf q}= {\bf q}'=(9,6,1,1,1).
$$

Notice that $q_i' >0$ for any $i$. Moreover, using Magma we found that
$$
\Gamma({\bf q}) \cong \left( {\mathbb Z}/18{\mathbb Z}\right)^3;
$$
$$
H_A \cong {\mathbb Z}/6{\mathbb Z} \times {\mathbb Z}/18{\mathbb Z};
$$
$$
\Gamma_A \cong {\mathbb Z}/3{\mathbb Z} \times {\mathbb Z}/18{\mathbb Z}.
$$

The Calabi-Yau $X_A$ has a singularity at $[-1,1,0,0,0]$, which is a singularity of $W{\mathbb P}^4(9,6,1,1,1)$. As explained in \cite{Karp}, this singularity can be blown-up so as to get a Calabi-Yau in the (toric) blow-up of $W{\mathbb P}^4(9,6,1,1,1)$.

{\it Example D}. When the group $H_A$ is trivial, it is possible to
write down an explicit birational inverse between from
$\overline{M}_A$ to $X_A$. We show it in one specific example. Let us
consider the polynomial
$$
F_A:=x_1^5+x_2^9x_3+x_3^9x_4+x_4^{10}+x_5^{2}
$$
where $A$ is the matrix $$
A:=\left(
\begin{array}{ccccc}
5 & 0 & 0 & 0 & 0 \\
0 & 9 & 1 & 0 & 0 \\
0 & 0 & 9 & 1 & 0 \\
0 & 0 & 0 & 10 & 0 \\
0 & 0 & 0 & 0 & 2
\end{array}
\right).
$$

The variety $X_A:=Z(F_A)$ is contained in $W{\mathbb P}^4(2,1,1,1,5)$. The matrix $B$ is given by
$$
B:=\left(
\begin{array}{ccccc}
162 & 0 & 0 & 0 & 0 \\
0 & 90 & -10 & 1 & 0 \\
0 & 0 & 90 & -9 & 0 \\
0 & 0 & 0 & 81 & 0 \\
0 & 0 & 0 & 0 & 405
\end{array}
\right),
$$
where $AB=BA=810I$. The map $q_A: W{\mathbb P}^4(2,1,1,1,5)
\rightarrow {\mathbb P}^5$
$$
(x_1: \ldots : x_5) \rightarrow (u_0: \ldots : u_5):= (x_1x_2\ldots
x_5: x^5_1: x_2^9x_3: x_3^9x_4:x_4^{10}x_2:x_5^{2})
$$
maps $X_A$ to the variety
$$
\overline{M}_A:=Z(u_1+\ldots+u_5,
-u^{810}+u_1^{162}u_2^{90}u_3^{80}u_4^{73}u_5^{405}) \subset
{\mathbb P}^5.
$$

An explicit birational inverse for $q_A$ is given by the following
map:

$$
\left\{
\begin{array}{l}
M^{162}x_1=u_1^{65}u_2^{54}u_3^{12}u_4^{15}u_5^{162}, \\
\\
M^{81}x_2=u_1^{-2}u_2^8u_3^{-11}u_4^{-8}u_5^{-5}, \\
\\
M^{81}x_3=u_1^{18}u_2^{19}u_3^{-1}u_4u_5^{45}, \\
\\
M^{81}x_4=u_2^9u_3^{-10}u_4^{-7}, \\
\\
M^{405}x_5=u_1^{81}u_2^{90}u_3^{-10}u_4u_5^{203},
\end{array}
\right.
$$
where $M=x_1^2x_2^3x_3x_4x_5^2$.

\section{A one-dimensional family}\label{sec4}

Let us consider the one-parameter family ${\mathcal X} \rightarrow {\mathbb P}^1_t$ of degree $d$ hypersurfaces
in ${\mathbb P}^{n-1}$ with ${\mathcal X}_t=X_{d,t}=Z(F_{d,t})$, where
\begin{equation}
\label{mirror}
F_{d,t}= \sum_i^n y_i^d - t \, y_1^{q_1'}\ldots y_n^{q'_n}.
\end{equation}

Clearly, this is a one-dimensional deformation of the Fermat variety $X_d$.  If we apply the map $\phi_A: X_{d,t} \subset
{\mathbb P}^{n-1} \rightarrow X_{A,t} \subset W {\mathbb P}(q_1, \ldots, q_n)$, the image of $X_{d,t}$ is given by
$
X_{A,t}=Z(F_{A,t}),
$ where
$$
F_{A,t} = F_A- t\, x_1x_2\ldots x_n.
$$

Under the composition $q_A \circ \phi_A$ the image
 of \eqref{mirror} is given by the equations
 \begin{equation}
 \label{uzero}
\sum_i u_i - t u_0=0, \qquad u^d_0 = u_1^{q_1'}\ldots u_n^{q_n'}.
\end{equation}

If $t \neq 0$, we solve for $u_0$ and get the equation
\begin{equation}
\label{mirr}
\left( \sum_i ^n u_i \right)^d = t^d u_1^{q_1'} \ldots u_n^{q_n'}.
\end{equation}

The group $\Gamma({\bf q'})$ acts on each $X_{A,t}$ since $\sum_i q'_i$ equals $d$. Denote by $\overline{\mathcal M} \rightarrow {\mathbb P}^1_t$ the family
given by \eqref{mirr}. By the universal property of the quotient there exists a map $\Psi$ between
${\mathcal X}/\Gamma({\bf q'})$
and ${\overline {\mathcal M}}$, which commutes with the projection map on ${\mathbb P}^1_t$.

\begin{proposition}
The map $\Psi$ yields a birational morphism from ${\mathcal X}/ \Gamma({\bf q}')$ to ${\overline{\mathcal M}}_t$.
\end{proposition}
\proof It suffices to compare the degree of the quotient map ${\mathcal X} \rightarrow {\mathcal X}/\Gamma({\bf q'})$, which is $\#\Gamma({\bf q'})$, with that
of the map ${\mathcal X} \rightarrow \overline{{\mathcal M}}$. Let
$$
(l_0:l_1 \ldots: l_n)=(\prod_k^ny_k^{q_k'}: y_1^d : \ldots y_n^d)
$$
be a generic point in $\overline{Im( q_A \circ \phi_A)}$. Thus, we have $y_j=\zeta^{k_j}\sqrt[d]{l_j}$, where $\zeta$ is a primitive
 $d$-th root of unity. Hence, we get
 $$
 \zeta^{\sum_j q_j'k_j}\sqrt[d]{l_1^{q_1'}\ldots l_n^{q_n'}}=l_0.
 $$

On the other hand, by the equation of $F_{A,t}$,  we must have
\begin{equation}
\label{set}
 \sum_j q_j'k_j \equiv 0 \, \, \hbox{mod} \, \,  d.
\end{equation}

Recall that $\sum_j q_j' \equiv 0$ mod $d$, so we can take the quotient of the set \eqref{set} of solutions
by $(1, 1, \ldots, 1)$. The degree of ${\mathcal X} \rightarrow {\overline {\mathcal M}}$ is thus equal to $\# \Gamma({\bf q'})$.

\endproof

\subsection{Some Birational Families}\label{sec8} Let us examine some Calabi-Yau varieties $X_A$ which have a one-dimensional versal deformation family,
so it may be described by the family in the section above.  In
Schimmrigk's list (see \cite{sch}), we found twelve entries with
$h^{2,1}=1$.  Since the deformation space is one dimensional, we
take into account one-dimensional families corresponding to these
entries of the list. The generic members of the families, which are denoted by
the same letter, have the same Euler characteristic. We have

$$
{\mathcal
A}_1(t)=x_1^8x_3+x_2x_3^7+x_2^7x_4+x_4^7+x_5^2-tx_1x_2x_3x_4x_5
\subset {\mathbb A}^1 \times W{\mathbb P}^4(75,84,86,98,343),
$$
$$
{\mathcal
A}_2(t)=x_1^8x_2+x_2^7x_3+x_3^7+x_4^4+x_4x_5^2-tx_1x_2x_3x_4x_5
\subset {\mathbb A}^1 \times W{\mathbb P}^4(43,48,56,98,147),
$$
$$
{\mathcal
A}_3(t)=x_2^8+x_1^7x_3+x_3^7+x_1x_4^4+x_4x_5^2-tx_1x_2x_3x_4x_5,
\subset {\mathbb A}^1 \times W{\mathbb P}^4(48,49,56,86,153)$$
$$
{\mathcal
A}_4(t)=x_1x_2^7+x_1^7x_3+x_3^7+x_2x_4^4+x_4x_5^2-tx_1x_2x_3x_4x_5\subset
{\mathbb A}^1 \times W{\mathbb P}^4(42,43,49,75,134),
$$

$$
{\mathcal
B}_1(t)=x_1^{10}x_2+x_2^9x_3+x_3^9+x_4^5+x_5^2-tx_1x_2x_3x_4x_5
\subset {\mathbb A}^1 \times W{\mathbb P}^4(73,80,90,162,405),
$$
$$
{\mathcal
B}_2(t)=x^9_1x_2+x_2^9+x_1x_3^8+x_3x_4^5+x_5^2-tx_1x_2x_3x_4x_5
\subset {\mathbb A}^1 \times W{\mathbb P}^4(64,72,73,115,324),
$$
$$
{\mathcal
B}_3(t)=x_1^9x_2+x_2^9+x_1x_3^5+x_4^5+x_3x_5^2-tx_1x_2x_3x_4x_5\subset
{\mathbb A}^1 \times W{\mathbb P}^4(40,45,73,81,166),
$$
$$
{\mathcal
B}_4(t)=x_1^9+x_2^8+x_2x_3^5+x_1x_4^5+x_4x_5^2-tx_1x_2x_3x_4x_5\subset
{\mathbb A}^1 \times W{\mathbb P}^4(40,45,63,64,148),
$$

$$
{\mathcal
C}_1(t)=x_1^6x_3+x_2x_3^5x_2^5+x_4+x_4^5+x_5^3-tx_1x_2x_3x_4x_5\subset
{\mathbb A}^1 \times W{\mathbb P}^4(52,60,63,75,125),
$$
$$
{\mathcal
C}_2(t)=x_2^6+x_1^5x_3+x_3^5+x_1x_4^4+x_4x_5^3-tx_1x_2x_3x_4x_5
\subset {\mathbb A}^1 \times W{\mathbb P}^4(48,50,60,63,79),
$$

$$
{\mathcal
D}_1(t)=x_1^5x_3+x_3^4x_4+x_2x_4^4+x_2^4x_5+x_5^4-tx_1x_2x_3x_4x_5\subset
{\mathbb A}^1 \times W{\mathbb P}^4(41,48,51,52,64),
$$
$$
{\mathcal
D}_2(t)=x_3^5+x_1^5x_4+x_2x_4^4+x_2^4x_5+x_5^4-tx_1x_2x_3x_4x_5\subset
{\mathbb A}^1 \times W{\mathbb P}^4(51,60,64,65,80).
$$

Let $V_5(t), V_6(t), V_8(t), V_{10}(t)$ be the four hypergeometric families on page 134
in \cite{mey}. For each of these families $V_j(t)$, $j=5,6,8,10$,
there is a group acting on the family such that the mirror $W_j(t)$
of $V_j(t)$ can be described as a resolution of the quotient
$V_j(t)/G$. The singular members have one orbit of ordinary nodes
under the action of $G$ and the resolution of the quotient is a
rigid Calabi-Yau threefold, i.e., $h^{2,1}=0$.

\begin{theorem}
The following birational equivalences hold:
$$
{\mathcal A_1}(t)  \approx {\mathcal A_2}(t) \approx {\mathcal
A_3}(t)\approx {\mathcal A_4}(t) \approx W_8(t),
$$
$$
{\mathcal B_1}(t)  \approx {\mathcal B_2}(t) \approx {\mathcal
B_3}(t)\approx {\mathcal B_4}(t) \approx W_{10}(t),
$$
$$
{\mathcal C_1}(t)  \approx {\mathcal C_2}(t) \approx W_6(t),
$$
$$
{\mathcal D_1}(t)  \approx {\mathcal D_2}(t) \approx W_5(t).
$$
\end{theorem}

\proof It is easy to check that the general member of the families
above is a singular Calabi-Yau in four weighted projective space.
The singular locus is a rational curve. For some values of $t$ there
are extra singularities that are ordinary nodes, namely:
\begin{center}
\begin{tabular}{cccc}
${\mathcal A}_i(t)$ & $t^8=2^{16}$ \vline & ${\mathcal B}_i(t)$ &
$t^{10}=800000$ \\ \hline ${\mathcal C}_i(t)$ & $t^6=3^62^4$ \vline
& ${\mathcal D}_i(t)$ & $t^5=5^5$ \\ \hline
\end{tabular}
\end{center}

Let us focus on the two families ${\mathcal D}_1(t)$ and ${\mathcal
D}_2(t)$. The proof for the other cases can be dealt with
analogously. Define the two families
$$
X_{d_1,t}=\left\{\sum_{i=1}^5 y^{320}_i - \prod_i
y_i^{64}=0\right\},
$$
$$
X_{d_2,t}=\left\{\sum_{i=1}^5 y^{1280}_i - \prod_i
y_i^{256}=0\right\}.
$$

Let $A_1$ and $A_2$ be the following matrices:
$$
A_1:=\left(
\begin{array}{ccccc}
5 & 0 & 1 & 0 & 0 \\
0 & 0 & 4 & 1 & 0 \\
0 & 1 & 0 & 4 & 0 \\
0 & 4 & 0 & 0 & 1 \\
0 & 0 & 0 & 0 & 4
\end{array}
\right), \quad A_2:=\left(
\begin{array}{ccccc}
0 & 0 & 5 & 0 & 0 \\
5 & 0 & 0 & 1 & 0 \\
0 & 1 & 0 & 4 & 0 \\
0 & 1 & 0 & 0 & 4 \\
0 & 0 & 0 & 0 & 4
\end{array}
\right).
$$

Using Magma we found that the groups $H_{A_1}$ and $H_{A_2}$ are trivial - this
happens for all $12$ families.

We thus have
$$
q_{A_1} \circ \phi_{A_1} : X_{d_1, t} \rightarrow \overline{M}_{A_
1,t}$$
$$
(y_1: y_2: \ldots: y_5) \rightarrow  \left(\prod_i y_i^{64}:
y_1^{320}: \ldots : y_5^{320}\right)
$$
and, similarly,
$$
q_{A_2} \circ \phi_{A_2} : X_{d_2, t} \rightarrow \overline{M}_{A_
2,t}
$$
$$
(y_1: y_2: \ldots: y_5) \rightarrow  \left(\prod_i y_i^{256}:
y_1^{1280}: \ldots : y_5^{1280}\right)
$$

It is easy to check that $\overline{M}_{A_1,t} \cong
\overline{M}_{A_2,t} \cong \left\{ \sum_i u_i -tu_0, \,
u_0^5=u_1\ldots u_5\right\} \approx W_5(t)$.

Since $H_{A_1}$ and $H_{A_2}$ are trivial, ${\mathcal D}_1(t)$ and
${\mathcal D}_1(t)$ are birational since they are both birational to
$W_5(t)$.

\endproof

\subsection{Picard-Fuchs equations}\label{sec9}  When $X_A$ is a Calabi-Yau hypersurface,
the Hodge number $h^{2,1}(X_A)$ gives the number of independent
parameters of deformations of $X_A$. There exists a system of
partial differential equations, the so called GKZ-hypergeometric
system (see \cite{GKZ}), which yield Picard-Fuchs equations for the
variation of periods along families with central fiber $X_A$. When
$h^{2,1}(X_A)=1$, the Picard-Fuchs equation can also be found via a
generalization of the Griffiths-Dwork method for hypersurfaces in
weighted projective space: see, for instance, \cite{mor}.

\section{${\overline M}_A$ and the mirror family of Calabi-Yau hypersurfaces in weighted projective space}\label{sec5}

Assume $X_A$ is a Calabi-Yau manifold (as defined in Section 2) in
weighted projective space $W{\mathbb P}^{n-1}(q_1, \ldots, q_n)$,
where $q_i | Q$ and $Q=\sum_j q_j$.  Batyrev's mirror construction
(see, for instance, \cite{ph}) depends only on the polytope $\Delta$
associated to the toric variety $W{\mathbb P}^{n-1}(q_1, \ldots,
q_n)$ and not on the matrix $A$. As explained, for instance in
\cite{ph}, the Calabi-Yau varieties in the mirror family ${\mathcal
W} \rightarrow {\mathbb P}^1_x$ of a general section of the
anticanonical bundle ${\mathcal O}(Q)$, with $Q=\sum_j q_j$,  can be
represented as compactifications of complete intersections of the
affine hypersurfaces in $({\mathbb C}^*)^{n}$ given by
\begin{equation}
\label{okkei}
t_1 + \ldots + t_n=1, \qquad t_1^{q_1}\hdots t_n^{q_n}=x.
\end{equation}

Let $W_x$ be the fiber over the point $x \in {\mathbb P}^1$. By comparing \eqref{uzero} and \eqref{okkei}, the following holds.

\begin{proposition}
\label{mirro}
The compactification of ${\mathcal W}_1$ is given by
the equations \eqref{uzero} that define the Shioda quotient ${\overline M}_{{}^tA,1}$ for any matrix $A$.
\end{proposition}

\proof Let $A$ be a matrix as in Section \ref{sec1}. If we start from the family $F_{{}^tA, t}$,
the equations \eqref{uzero} become:
\begin{equation}
\sum_i u_i - t u_0=0, \qquad u^d_0 = u_1^{q_1}\ldots u_n^{q_n}.
\end{equation}

Since $\sum_j q_j=d$, the claim follows.
\endproof

\medskip

{\bf Acknowledgements.} The author wishes to express his gratitude
to Bert van Geemen for helpful remarks and useful suggestions.

 \end{document}